\documentclass[11pt,a4paper]{article}
\usepackage[english]{babel}
\usepackage{amsmath}
\usepackage{amsfonts}
\usepackage{amsthm}
\usepackage{graphicx}
\usepackage{makeidx}

\title{ Instability of Hopf vector fields on Lorentzian Berger spheres}
\author{
 Ana Hurtado \thanks{Partially supported
by DGI (Spain) and FEDER Project MTM 2004-06015-C02-01  and by
Generalitat Valenciana Grant ACOMP06/166.}\\ [4pt] Departament de
Matem\`atiques, Universitat Jaume I \\ E-12071 Castell\'o, Spain.\\
e-mail: ahurtado@mat.uji.es}
\date{}

\theoremstyle{plain} 
\newtheorem{thm}{Theorem}[section]
\newtheorem{cor}[thm]{Corollary}
\newtheorem{lem}[thm]{Lemma}
\newtheorem{prop}[thm]{Proposition}

\newtheorem{defn}[thm]{Definition}
\theoremstyle{definition}

\theoremstyle{remark}

\newtheorem{nota}[thm]{Remark}
\begin{document}
\maketitle

\begin{abstract}
\parindent0pt\noindent
In this work, we study the stability of Hopf vector fields on
Lorentzian Berger spheres as critical points of the energy, the
volume and the generalized energy. In order to do so, we construct
a family of vector fields using the simultaneous eigenfunctions of
the Laplacian and of the vertical Laplacian of the sphere. The
Hessians of the functionals are negative when they act on these
particular vector fields and then Hopf vector fields are unstable.
Moreover, we use this technique to study some of the open problems
in the Riemannian case.



\end{abstract}

\section{Introduction}
A smooth vector field $V$ on a Riemannian manifold $(M,g)$ can be
seen as a map into its tangent bundle endowed with the Sasaki
metric, $g^S$, defined by $g$. The volume of $V$ is the volume of
$V(M)$ considered as a submanifold of $(TM, g^S)$ . Analogously,
we can define the energy of $V$ as the energy of the map
$V:(M,g)\longrightarrow (TM, g^S)$ and if $\tilde g$ is another
metric on $M$, we define the generalized energy $E_{\tilde g}$ as
the energy of $V:(M, \tilde g)\to (TM,g^S)$. These energies were
introduced in \cite{GM} to study the relationship between the
energy and the volume of vector fields. In particular, if we take
either $\tilde g=g$ or $\tilde g= V^* g^S$, the generalized energy
turns out to be, up to constant factors, the energy and the volume
of the vector field respectively.\par

On a compact manifold $M$, the critical points of all these
functionals should be parallel with respect to the Levi-Civita
connection defined by $g$, so it is usual to restrict the
functionals to the submanifold of unit vector fields. Obviously,
if $M$ admits unit parallel vector fields, they are the absolute
minimizers.\par

The geometrically simplest manifolds admitting unit vector fields
but not parallel ones are odd-dimensional spheres. Hopf vector
fields defined as those tangent to the fibres of the Hopf
fibration $\pi: S^{2m+1}\longrightarrow \mathbb{C}P^m$ are very
special unit vector fields. When both manifolds are endowed with
their usual metrics, this map is a Riemannian submersion with
totally geodesic fibres whose tangent space is generated by the
unit vector field $V= JN$, where $N$ is the unit normal to the
sphere and $J$ is the usual complex structure of $\mathbb
R^{2m+2}$.\par

In \cite{G-Z}, Gluck and Ziller showed that Hopf vector fields on
the $3$-dimensional round spheres are the absolute minimizers of
the volume and the analogous result for the energy was shown by
Brito in \cite{Br}. For higher dimension, they are unstable
critical points of the energy (\cite{GMLF2}, \cite{Wo} and
\cite{Wo2}).\par

All these results are independent of the radius of the sphere, but
as concerns the stability as critical points of the volume,
Borrelli and Gil-Medrano showed in \cite{BoGM} that for each $m>1$
there exists a critical value of the radius, such that, Hopf
vector fields are stable critical points of the volume if and only
if the radius is lower than or equal to this critical radius. By
stable we mean that the Hessian of the functional is positive
semi-definite.\par

In order to understand better these phenomena, in \cite{GMH}
Gil-Medrano and the author studied the behaviour of the Hopf
vector field with respect to the volume and the energy when the
metric considered on the sphere is the canonical variation of the
Riemannian submersion given by the Hopf fibration. The metrics so
constructed are known as Berger metrics, they consist in a
$1$-parameter variation $g_\mu$ for $\mu\neq 0$. When $\mu>0$, the
new metric is Riemannian and if $\mu<0$, the metric is Lorentzian
and $V^\mu=1/\sqrt{-\mu}\, V$ is timelike. Moreover, they also
studied the subset of $\mathbb{R}^+\times \mathbb{R}^+$ of pairs
$(\mu,\lambda)$ such that the vector field $V^\mu$ is stable as a
critical point of the generalized energy $E_{g_\lambda}$ on the
spheres of dimension greater than three. The dimension three was
studied in \cite{AH}.\par

In Riemannian Berger spheres, the problem of determining the
behaviour of Hopf vector fields is completely solved for the
energy and the volume, but as concerns the generalized energy
$E_{g_\lambda}$, there exist values of $\lambda$ and $\mu$ for
which the stability of Hopf vector fields is still an open
problem. For Lorentzian Berger spheres, the technique used to show
stability in the Riemannian case does not allow us to conclude the
stability in any case and only a partial result concerning  the
instability is shown in \cite{GMH}. These instability results, as
in the Riemannian case, have been obtained computing the Hessian
in the direction of the vector fields $A_a=a-\langle a, V\rangle
V-\langle a, N\rangle N$ for all $a\in \mathbb{R}^{2m+2}$, $a\neq
0$. These vector fields can be seen as the projection onto
$V^\bot$ of the gradient of an eigenfunction associated to the
first eigenvalue of the Laplacian of the sphere.\par

In this work, we construct new directions using the simultaneous
eigenfunctions of the Laplacian  and of the vertical Laplacian
$\Delta_v (f)=-V(V(f))$ of the sphere. More precisely, we consider
vector fields $C_{2s}={\rm grad}^\mu f_{2s}-\varepsilon_\mu
V^\mu(f_{2s})f_{2s}$, where $f_{2s}$ is a polynomial of degree
$2s$ in $\mathbb{R}^{2m+2}$ such that its restriction to the
sphere is a simultaneous eigenfunction of the Laplacian and of the
vertical Laplacian. Here $\varepsilon_\mu=\mu/|\mu|$. These vector
fields verify that
$\nabla^\mu_{V^\mu}C_{2s}=(\mu-2s)/\sqrt{|\mu|}\ JC_{2s}$ and they
allow us to prove in Section $3$ that {\it on Lorentzian Berger
spheres, the Hopf vector fields $V^\mu$ are unstable critical
points of the energy, the volume and the generalized energy
$E_{g_\lambda}$ for all $\lambda <0$}. The eigenfunctions of
$\Delta_v$ have been also used to study, for example, the harmonic
index and nullity of the Hopf map (see \cite{L}).\par

In Section $4$, we use the ideas introduced in the previous
section to complete the results in \cite{AH} and then we solve
completely the problem of determining the stability of Hopf vector
fields with respect to the generalized energy $E_{g_\lambda}$ in
the Riemannian Berger $3$-sphere. In particular, we prove that
{\it if $\lambda
> (\mu-3)^2/(\mu-2)$ and $\mu>2$, or if $\lambda
> \mu-4$ and $\mu>4$, then $V^\mu$ is an unstable critical point
of $E_{g_\lambda}$}. Again, we need to consider vector fields more
complicated that the vectors fields $A_a$. So, the simultaneous
eigenfunctions of the Laplacian and of the vertical Laplacian play
an important role in the resolution of these problems in the
sphere. For spheres of upper dimension, we can use the vector
fields $C_{2s}$ to improve the results in \cite{GMH} concerning
the generalized energy, but it is not sufficient to solve
completely the problem.

\section{Preliminaries}

Given a Riemannian manifold $(M,g)$, the Sasaki metric $g^S$ on
the tangent bundle $TM$ is defined, using $g$ and its Levi-Civita
connection $\nabla$,  as follows:
$$
g^S(\zeta_1,\zeta_2)=g(\pi_* \circ \zeta_1,\pi_* \circ
\zeta_2)+g(\kappa\circ \zeta_1, \kappa\circ\zeta_2),
$$
where $\pi: TM\to M$ is the projection and $\kappa$ is the
connection map of $\nabla$. We will consider also its restriction
to the tangent sphere bundle, obtaining the Riemannian manifold
$(T^1M, g^S).$\par As in \cite{GM}, for each metric $\tilde g$ on
$M$ we can define the generalized energy of the vector field $V$,
denoted $E_{\tilde g}(V)$, as the energy of the map $V: (M,\tilde
g)\to (TM, g^S)$ that is given by
$$
E_{\tilde g}(V) =\frac 12\int_M \operatorname{tr} L_{(\tilde
g,V)}\,{\rm dv}_{\tilde g} ,$$ where $L_{(\tilde g,V)}$ is the
endomorphism determined by $V^*g^S (X,Y) = \tilde g(L_{(\tilde
g,V)}(X),Y)$. This energy can also be written as
\begin{equation}\label{energiag}
E_{\tilde g}(V) = \frac12 \int_M \sqrt{ \operatorname{det}
P_{\tilde g}}\, \operatorname{tr} (P_{\tilde g}^{-1}\circ L_V)\,
{\rm dv}_g,
\end{equation}
where $P_{\tilde g}$  and $L_V$ are defined by $\tilde g(X,Y) =
g(P_{\tilde g}(X),Y)$ and $V^*g^S (X,Y) = g(L_V(X),Y)$,
respectively. By the definition of the Sasaki metric, $L_V =
\operatorname{Id} + (\nabla V)^t\circ\nabla V$. In particular, for
$\tilde g=g$
\begin{equation}\label{energia}
E_g(V)=\frac12\int_M  \operatorname{tr} L_V\,{\rm dv}_g=\frac n2
vol(M,g)+\frac12\int_M \Vert \nabla V \Vert^2\,{\rm dv}_g.
\end{equation}
This functional is known as the energy and will be represented by
$E$. Its relevant part, $B(V)=\frac12\int_M \Vert \nabla V \Vert^2
\,{\rm dv}_g$, is known as the total bending of $V$ and its
restriction to unit vector fields has been widely studied by
Wiegmink in \cite{Wie1}, (see also \cite{Wo}).\par

On the other hand, the volume of a vector field $V$ is defined as
the $n$-dimensional volume of the submanifold $V(M)$ of $(TM,
g^S)$. It is given by

\begin{equation}\label{volumen}
F(V)=\int_M \sqrt{\operatorname{det} L_V}\,{\rm dv}_g.
\end{equation}
Since for $\tilde g =V^*g^S$ we have $P_{\tilde g} =  L_V$, then
(\ref{energiag}) and (\ref{volumen}) give
$$
F(V)=\frac 2n E_{V^*g^S}(V).
$$
The first variation of the generalized energy has been computed in
\cite{GM}. It has been also shown there that $V$ is a critical
point of $F$ if and only if $V$ is a critical point of
$E_{V^*g^S}$ and that, on a compact $M$, a critical vector field
of any of these generalized energies should be parallel. So, it is
usual to restrict these functionals to the submanifold of unit
vector fields.\par

The following proposition shown in \cite{GM} generalizes the
characterization of critical points of the total bending in
\cite{Wie1} and of the volume in \cite{GMLF2}.\par
\begin{prop} \label{var1}  Let $(M, g)$ be a
Riemannian manifold, a unit vector field $V$ is a critical point
of $E_{\tilde g}$ if and only if
$$
\omega_{(V, \tilde g)}\, (V^\bot)=\{ 0 \},
$$
with $\ \omega_{(V, \tilde g)}=C_1^1 \nabla K_{(V, \tilde g)}\ $
and $\ K_{(V, \tilde g)}=\sqrt{\operatorname{det} P_{\tilde g}}\
P_{\tilde g}^{-1}\circ (\nabla V)^t.$
\end{prop}

\begin{nota} For a $(1,1)$-tensor field $K$,  if $\{E_i\}$ is a
$g$-orthonormal local frame,
$$C_1^1\nabla K (X) = \sum_i g((\nabla_{E_i}K)X,E_i).$$
\end{nota}
Moreover, in \cite{GMLF} it was proved that a unit vector field is
a critical point of $F$ if and only if it defines a minimal
immersion in $(T^1 M,g^S)$.
\begin{thm}[\cite{GMLF2}]\label{var2} Let $V$ be a unit
vector field on the Riemannian manifold $(M,g)$.\par
\smallskip

\noindent a)  If $V$ is a critical point of $E_{\tilde g}$, the
Hessian of $E_{\tilde g}$ at $V$ acting on $A\in V^\bot$ is given
by
$$
(Hess E_{\tilde g})_V(A) =\int_M \Vert A \Vert^2 \omega_{(V,
\tilde g)}\,(V)\,{\rm dv}_g +\int_M
\sqrt{\operatorname{det}P_{\tilde g}}\ \operatorname{tr} \Big(
P_{\tilde g}^{-1}\circ(\nabla A)^t \circ\nabla A \Big)\,{\rm dv}_g
.
$$

\noindent b) If $V$ is a critical point of the energy, the Hessian
of $E$ at $V$ acting on $A\in V^\bot$ is given by
$$
(Hess E)_V(A)=\int_M \Vert A \Vert ^2 \omega_{(V,g)}(V)\,{\rm
dv}_g+ \int_M \Vert \nabla A\Vert ^2\,{\rm dv}_g.
$$

\noindent c) Let $V$ be a unit vector field defining a minimal
immersion, the Hessian of $F$ at $V$ acting on $A\in V^\bot$  is
given by
\begin{eqnarray*}
(Hess F)_V(A)\ =& &\!\!\int_M \Vert A \Vert^2 \omega_V(V)\,{\rm
dv}_g
 + \int_M  \frac 2 {\sqrt{\operatorname{det} L_V}}\ \sigma_2 (K_V
\circ \nabla A)\,{\rm dv}_g\\  &-&\int_M  \operatorname{tr} \Big(
L_V ^{-1} \circ (\nabla A)^t \circ  \nabla V
\circ  K_V \circ  \nabla A \Big)\,{\rm dv}_g \\
 &+&  \int_M \sqrt{\operatorname{det} L_V}\ \operatorname{tr} \Big(
L_V ^{-1}\circ  (\nabla A)^t\circ  \nabla A \Big)\,{\rm dv}_g,
\end{eqnarray*}
where $\sigma_2$ is the second elementary symmetric polynomial
function. In particular, $\sigma_2 (K_V\circ \nabla A) =
(\operatorname{tr}(K_V\circ \nabla A))^2 -
\operatorname{tr}(K_V\circ \nabla A)^2$.
\end{thm}
\smallskip

The generalized energy can be defined for any $g$ and $\tilde{g}$
semi-Riemannian metrics on the manifold $M$. In particular, in a
Lorentzian manifold, the energy is defined for all vector fields.
Nevertheless, the volume of a reference frame (unit timelike
vector field) $V$ is not always defined, since the $2$-covariant
field $V^* g^S$ can be degenerated. Due to this, we study the
volume restricted to unit timelike vector fields for which $V^*
g^S$ is a Lorentzian metric on $M$. We will denote this set of
vector fields by $\Gamma^-(T^{-1}M)$ and it is an open subset of
the set of smooth references frames. If $V$ belongs to
$\Gamma^-(T^{-1}M)$, then ${\rm det}\, L_V > 0$ and the volume is
well defined.\par

The condition for a reference frame to be a critical point of the
generalized energy on a Lorentzian manifold is the same condition
that the one given by Proposition \ref{var1} for Riemannian
metrics. If we compute the second variation, we obtain the
following
\begin{prop}\label{varLor}  Let $V$ be a unit timelike vector field
on a compact Lorentzian manifold $(M,g)$.\par
\begin{enumerate}
\item If $V$ is a critical point of the generalized energy
$E_{\tilde{g}}$ and $A\in V^\bot$, then
\begin{equation}
 \label{hesenerlor}(Hess E_{\tilde{g}})_V(X)=-\int_M \|X\|^2\omega_{(V,\tilde{g})}(V)\,{\rm
 dv}_g+\int_M {\rm tr}(L_{\tilde{g}}\circ(\nabla X)^t\circ
 \nabla X)\,{\rm
 dv}_g.
 \end{equation}

 \item \cite{H} If $V$ is a critical
point of the energy, the Hessian of $E$ at $V$ acting on $A\in
V^\bot$ is given by
\begin{displaymath}
(Hess E)_V(A)=-\int_M \Vert A \Vert ^2 \omega_{(V,g)}(V)\, {\rm
dv}+ \int_M \Vert \nabla A\Vert ^2\, {\rm dv}.
\end{displaymath}

\item \cite{H} For a unit timelike vector field $V\in
\Gamma^{-}(T^{-1}M)$ defining a minimal immersion, the Hessian of
$F$ at $V$ acting on $A\in V^\bot$  is given by
\begin{eqnarray*}
(Hess F)_V(A)\ =&-&\int_M \Vert A \Vert^2 \omega_V(V)\, {\rm dv}
 + \int_M  \frac 2 {\sqrt{\operatorname{det} L_V}}\ \sigma_2 (K_V
\circ \nabla A)\,{\rm dv}\\  &-&\int_M  \operatorname{tr} \Big(
L_V ^{-1} \circ (\nabla A)^t \circ  \nabla V
\circ  K_V \circ  \nabla A \Big)\, {\rm dv} \\
 &+&  \int_M \sqrt{\operatorname{det} L_V}\ \operatorname{tr} \Big(
L_V ^{-1}\circ  (\nabla A)^t\circ  \nabla A \Big)\, {\rm dv}.
\end{eqnarray*}
\end{enumerate}
\end{prop}

The expression of the Hessian of the generalized energy given by
\eqref{hesenerlor} is obtained by straightforward computation in a
similar way that in the Riemannian case, so we have omitted the
details.

\begin{nota}\label{nota loren}Let us point out that if we compare the above
expressions of the Hessian with those obtained for Riemannian
metrics, the only difference is the minus sign of the first term
of the expression of the Hessian.
\end{nota}

\bigskip

Hopf vector fields on odd-dimensional spheres are tangent to the
fibres of the Hopf fibration $\pi:(S^{2m+1},g)\to (\mathbb
C^m,\overline g)$, where $g$ is the usual metric of curvature $1$
and $\overline g$ is the Fubini-Study metric with sectional
curvatures between $1$ and $4$. This map is a Riemannian
submersion with totally geodesic fibres whose tangent space is
generated by the unit vector field $V= JN$, where $N$ is the unit
normal to the sphere and $J$ is the usual complex structure of
$\mathbb R^{2m+2}$; in other words,  $V(p) = i p$.\par

In $S^{2m+1}$ we can consider the canonical variation $g_\mu$,
with $\mu\neq 0$, of the usual metric $g$,

\begin{equation}\label{metBerger}
g_\mu\vert_{V^\bot} = g\vert_{V^\bot},\qquad g_\mu(V,V) = \mu g
(V,V),\qquad g_\mu(V,V^\bot) = 0.
\end{equation}
 When $\mu>0$
the new metric is Riemannian and if $\mu<0$ the metric is
Lorentzian and $V$ is timelike.\par
 For all $\mu\not =0$, the map
$\pi:(S^{2m+1},g_\mu)\to (\mathbb C^m,\overline g)$ is  a
semi-Riemannian submersion with totally geodesic fibres.
$(S^3,g_\mu)$, with $\mu>0$, is known as a Berger sphere. We will
use the same name for all dimension and we will call $V^\mu =
\frac{1}{\sqrt{\vert\mu\vert}}V$ the Hopf vector field. It is a
unit Killing vector field with geodesic flow.\par

We denote by $\bar\nabla$ the Levi-Civita connection on $\Bbb
R^{2m+2}.$ The Levi-Civita connection $\nabla$ on $(S^{2m+1}, g) $
is $\nabla_XY=\bar\nabla_XY-<\bar\nabla_XY,N>N$ and
$\bar\nabla_XV=J\bar\nabla_XN= JX.$ Therefore $\nabla_V V=0$
and if $<X,V>=0$ then $\nabla_XV= JX.$\\

Using Koszul formula, one obtains the relation of $\nabla^\mu$,
the Levi-Civita connection of the metric $g_\mu$, with $\nabla$
\begin{equation}\label{dcBerger}
\nabla^\mu_VX = \nabla_VX +(\mu -1) \nabla_XV ,\qquad
\nabla^\mu_XV =  \mu \nabla_XV ,\qquad \nabla^\mu_XY = \nabla_XY,
\end{equation}
for all $X,Y$ in $V^\bot$.\par

It has been shown in \cite{GMH} that,
\begin{prop}\label{Hopfharmonic} For all $\mu, \lambda\not=0$, the
map $V^\mu: (S^{2m+1},g_\lambda)\to (T^1(S^{2m+1}), g_\mu^S)$ is
harmonic.
\end{prop}

Since $(V^\mu)^*g_\mu^S = (1+\vert\mu\vert)g_\lambda$ where
$\lambda = \mu/(1+\vert\mu\vert)$, as a consequence of the
Proposition above, we have the following
\begin{cor}[\cite{GMH}] For all $\mu\not=0$, the Hopf vector field
$V^\mu$ is a critical point of the generalized energy
$E_{g_\lambda}$, for all $\lambda\not=0$,  and it defines a
minimal immersion.
\end{cor}
\begin{nota}
When $\mu<0$, $V^\mu$ induces on the sphere a Lorentzian metric
$(V^\mu)^* g_\mu^S$ and the Hopf vector field is a critical point
of the volume restricted to the set of unit timelike vector fields
verifying this condition.
\end{nota}

The second variation of the energy and the volume at Hopf vector
fields on Berger spheres has been computed in \cite{GMH}. The
expression of the Hessian of the generalized energy
$E_{g_\lambda}$ is also computed in \cite{GMH} for Riemannian
Berger spheres. In a similar way, by straightforward computation,
we can obtain the second variation of the generalized energy
$E_{g_\lambda}$ in the Lorentzian case.\par

\begin{prop}\label{HesHopf}
Let $V^{\mu}$ be the Hopf unit vector field on
$(S^{2m+1},g_{\mu})$. For each vector field $A$ orthogonal to
$V^{\mu}$ we have
\begin{eqnarray*}
&&\textrm{a) }
    (Hess E_{g_{\lambda}})_{V^{\mu}}(A)=\int_{S^{2m+1}}
    \Big(-2m\varepsilon_\mu\sqrt{|\lambda\mu|}\ \|A\|^2+
    \sqrt{|\lambda/\mu|}\ \|\nabla^{\mu}A\|^2\\
    &&\hspace{1.6in}+\ (\varepsilon_\lambda\sqrt{|\mu/\lambda|}
  -\varepsilon_\mu\sqrt{|\lambda/\mu|}\ )\ \|\nabla^{\mu}_{V^{\mu}}A\|^2\Big)\,{\rm dv}_{\mu}.\\
&&\textrm{b) }
                (Hess E)_{V^{\mu}}(A)=\int_{S^{2m+1}}
    \Big(-2m\mu\|A\|^2+\|\nabla^{\mu}A\|^2\Big)\,{\rm dv}_{\mu}.\\
   &&\textrm{c) }
                (Hess F)_{V^{\mu}}(A)=(1+|\mu|)^{m-2}\int_{S^{2m+1}}
    \Big(\|\nabla^{\mu}A\|^2+ \ \mu\|\nabla^{\mu}_{V^{\mu}}A+\varepsilon_\mu\sqrt{|\mu|} JA\|^2\\
    && \hspace{1.6in}+\mu(-2m-2m|\mu|+2\varepsilon_\mu+2\varepsilon_\mu(m-\mu))\|A\|^2\Big)\,{\rm dv}_{\mu}.
\end{eqnarray*}
Where $\varepsilon_\mu=\mu/|\mu|$ and
$\varepsilon_\lambda=\lambda/|\lambda|$.
\end{prop}

\bigskip
Finally, let us recall some results concerning the vertical
Laplacian.\par

Let $\pi:(M,g)\longrightarrow (N,h)$ be a Riemannian submersion
and let $\Delta$ the Laplacian of $(M,g)$.
\begin{defn}[\cite{BBB}] The vertical Laplacian $\Delta_v$ of $(M,g)$ is the differential
operator given by
\begin{displaymath}
(\Delta_v f)(x)=(\Delta_{F_x}(f_{\vert_{F_x}}))(x),
\end{displaymath}
where $F_x=\pi^{-1}(\pi(x))$ is the fibre of $\pi$ passing through
 $x$ and $\Delta_{F_x}$ is the Laplacian of the induced metric by $M$ on $F_x$.\par
 The difference operator
$\Delta_h=\Delta-\Delta_v$ is called the horizontal Laplacian.
\end{defn}
B\'{e}rard-Bergery and Bourguignon, showed in \cite{BBB} that if
$\pi$ is a Riemannian submersion with totally geodesic fibres,
then $\Delta$ and $\Delta_v$ commute. So, when $M$ is compact and
connected, $L^2(M)$ admits a Hilbert basis consisting of
simultaneous eigenfunctions of both operators.

On $(S^{2m+1},g)$, it is known (see \cite{BEGM}) that the
eigenvalues of the Laplacian are $\lambda_k=k(k+2m)$ with
$k=0,1,2,\ldots$. Moreover, the eigenvalues of the vertical
Laplacian $\Delta_v$ are $\phi_l=l^2$ with $l=0,1,2,\ldots$. Then,
as can be seen in \cite{BBB} and \cite{HP}, using that the
Laplacian of the metrics $g_\mu$ is,
$\Delta^\mu=\mu^{-1}\,\Delta_v + \Delta_h$, the eigenvalues of
$\nabla^\mu$ are of the type

\begin{eqnarray}
\label{valores
propios}\lambda^\mu_{k,l}&=&(\lambda_k-\phi_l)+\frac{1}{\mu}\phi_l=k(k+2m)-l^2
+\frac{1}{\mu}l^2,\qquad k\geq l.
\end{eqnarray}
In the above expression not all values of $k$ and $l$ are
possible.\par Besides, Tanno showed that,
\begin{lem}[\cite{Tanno}]On $S^{2m+1}$, for each eigenvalue $\lambda_k$ of
$\Delta$, the space of eigenfunctions $\mathcal{P}^k$ admits an
orthogonal decomposition
\begin{displaymath}
\mathcal{P}^k=\mathcal{P}^k_k+\mathcal{P}^k_{k-2}+\cdots+\mathcal{P}^k_{k-2[k/2]},
\end{displaymath}
where $[k/2]$ is the integer part of $k/2$, and for $f \in
\mathcal{P}^k_{k-2p}$
\begin{displaymath}
V(V(f))=-(k-2p)^2f,\qquad 0\leq p\leq [k/2].
\end{displaymath}
Some of the $\mathcal{P}^k_{l}$ could be trivial.
\end{lem}

Since $\Delta_v(f)=-V(V(f))$, the problem of determining which
spaces $\mathcal{P}^k_l$ are not trivial is related to that of
determining the permitted combinations of $k$ and $l$ in
(\ref{valores propios}).\par In the following sections, we will
use that $\mathcal{P}^k_k\neq \{0\}$ for all $k$ (see
\cite{Tanno}), that is to say, for all $k>0$, $
\lambda^\mu_{k,k}=k(2m+\frac{1}{\mu}k)$, is an eigenvalue of
$\Delta^\mu$.

\section{Lorentzian Berger spheres }

The instability results in Riemannian Berger spheres have been
obtained by computing the Hessians when they act on the vector
fields $A_a=a-\langle a, V\rangle V-\langle a, N\rangle N$ for all
$a\in \mathbb{R}^{2m+2}$, $a\neq 0$. These vector fields can be
seen as the projection onto $V^\bot$ of the gradient of an
eigenfunction associated to the first eigenvalue of the Laplacian.
For Lorentzian Berger spheres, if we compute the Hessians in the
direction of these particular vector fields, we obtain
\begin{lem}
Let $V^{\mu}$ be the Hopf unit vector field on
$(S^{2m+1},g_{\mu})$, with $\mu<0$. For each $a \in
\mathbb{R}^{2m+2}$, $a\not = 0$ we have:
\begin{eqnarray*}
&&\textrm{a) }
                (Hess E)_{V^{\mu}}(A_a)=\frac{\sqrt{-\mu}m}{m+1}|a|^2
\Big((1-2m)\mu+2+\frac{(\mu-1)^2}{\mu}\Big){\rm vol}(S^{2m+1}).\\
&& \textrm{b) }
                (Hess
F)_{V^{\mu}}(A_a)=(1-\mu)^{m-2}\frac{\sqrt{-\mu}m}{m+1}|a|^2
f(m,\mu) {\rm vol}(S^{2m+1}),
\end{eqnarray*}
where $f(m,\mu) =
\Big((2m-1)\mu^2+(1-4m)\mu+2+(1-\mu)\frac{(\mu-1)^2}{\mu}\Big).$
\end{lem}
As a consequence,
\begin{prop}[\cite{GMH}]\label{InestL}
Let $V^{\mu}$ be the unit Hopf vector field on
$(S^{2m+1},g_{\mu})$, with $\mu<0$.
\begin{itemize} \item[a) ]If $(2m-2)\mu^2 <1$, then $V^\mu$ is
unstable for the energy. \item[b) ] If
$(2-2m)\mu^3+(4m-4)\mu^2+\mu <1$, then $V^\mu$ is unstable for the
volume. 
\end{itemize}
 In particular, on $(S^3, g_\mu)$ the Hopf vector field is unstable for the energy and the volume for all values of $\mu<0$.
\end{prop}

The alternative expressions of the Hessians used to show stability
in the Riemannian case (see \cite{GMH}) can be extended to include
negative values of $\mu$, but they do not allow us to conclude the
stability of Hopf vector fields in any case. In fact, we are going
to prove that they are always unstable. Moreover, we will study
the behaviour of Hopf vector fields with respect to the
generalized energy $E_{g_\lambda}$.\par

In order to do so, we are going to consider new directions
obtained from functions that are simultaneous eigenfunctions of
the Laplacian and of the vertical Laplacian of the sphere.\par

Let $f$ be a harmonic and homogeneous polynomial of degree $s$ in
$\mathbb{R}^{2m+2}$, then the restriction of $f$ to the sphere,
denoted also by $f$ for simplicity, is an eigenfunction associated
to the eigenvalue $s\,(2m+s)$ of the Laplacian of the sphere with
the usual metric. Moreover, we take $f$ verifying that
\begin{equation}
\label{JHess}\overline{Hess}f(u,Jv)=\overline{Hess}f(Ju,v)
\end{equation}
for all $u$, $v$ vector fields in $\mathbb{R}^{2m+2}$, where
$\overline{Hess}$ represents the Hessian in $\mathbb{R}^{2m+2}$.
In the sequel, we will denote with a bar the geometrical operators
related to the Euclidean space $\mathbb{R}^{2m+2}$.\par We will
use condition (\ref{JHess}) to assure that if $\{N,E_i,V,JE_i\}$
is a $J$-orthonormal local frame in $\mathbb{R}^{2m+2}$ then
$$\overline{Hess}f(E_i,E_i)+\overline{Hess}f(JE_i,JE_i)=0,\quad
\forall\ 1\leq i\leq m,$$ and
$$\overline{Hess}f(V,V)+\overline{Hess}f(N,N)=0.$$

\begin{prop}Let $f$ be a harmonic and homogeneous polynomial of degree $s$ in
$\mathbb{R}^{2m+2}$ satisfying (\ref{JHess}), then
\begin{itemize}
\item[a) ] If $\Delta^\mu$ denotes the Laplacian  of Berger
spheres and $\Delta^\mu_v$ the vertical Laplacian,
\begin{displaymath}
\Delta^\mu f=(2ms+\frac{s^2}{\mu})f\quad \textrm{ and }\quad
\Delta^\mu_v(f)=\frac{s^2}{\mu}f.
\end{displaymath}
In other words, $f \in \mathcal{P}^s_s$. \item[b) ] If $C={\rm
grad}^\mu f-\varepsilon_\mu V^\mu(f)V^\mu$, then
$\nabla^\mu_{V^\mu}C=(\mu-s)/\sqrt{|\mu|}\ JC$.
\end{itemize}
\begin{proof}
 If $u,v \in (V^\mu)^\bot$ then
\begin{eqnarray*}
    (Hess)^{\mu}f(V^{\mu},u)&=&V^{\mu}(u(f))-(\nabla^{\mu}_{V^{\mu}}u)f=\overline{Hess}f(V^{\mu},u)
    +\frac{1-\mu}{\sqrt{|\mu|}}Ju(f),\\
    (Hess)^{\mu}f(u,v)&=&u(v(f))-(\nabla^{\mu}_u
    v)f=\overline{Hess}f(u,v)-\langle u,v\rangle N(f),\\
    (Hess)^{\mu}f(V^{\mu},V^{\mu})&=&V^{\mu}(V^{\mu}(f))=\overline{Hess}f(V^\mu,V^\mu)
    -\frac{1}{|\mu|}N(f).
\end{eqnarray*}
Moreover, using (\ref{JHess}) and the fact that $N(f)=s\,f$, we
have that
\begin{displaymath}
\overline{Hess}f(V^\mu,V^\mu)=-\frac{\overline{Hess}f(N,N)}{|\mu|}=\frac{N(f)-N(N(f))}{|\mu|}=\frac{s(1-s)}{|\mu|}f.
\end{displaymath}
Then, since $f$ is a harmonic polynomial,
\begin{eqnarray*}
    \Delta^{\mu}(f)&=&-{\rm tr} (Hess )^{\mu}f=-\sum_{i=1}^{2m}(Hess
    )^{\mu}f(E_i,E_i)-\varepsilon_\mu \,(Hess
    )^{\mu}f(V^\mu,V^\mu)\\
    &=&-\sum_{i=1}^{2m}\overline{Hess
    }f(E_i,E_i)+2m N(f)+\frac{s^2}{\mu}f\\
    &=&(2ms+\frac{s^2}{\mu})f,
\end{eqnarray*}
and
$$\Delta^\mu_v(f)=-\varepsilon_\mu V^\mu(V^\mu(f))=-\varepsilon_\mu(\overline{Hess}f(V^\mu,V^\mu)-\frac{N(f)}{|\mu|})
=\frac{s^2}{\mu}f.$$

To show b), since $\nabla^\mu C=g_\mu^{-1}(Hess)^\mu
f-\varepsilon_\mu\,\nabla^\mu(V^\mu(f)V^\mu)$,
\begin{eqnarray*}
g_{\mu}(\nabla^{\mu}_{V^{\mu}}C,E_j)&=&(Hess )^{\mu}f(V^{\mu},E_j)
    -\varepsilon_{\mu}\,g_{\mu}(\nabla^{\mu}_{V^{\mu}}(V^{\mu}(f)V^{\mu}),E_j)\\
    &=&\overline{Hess}f(V^{\mu},E_j)+\frac{1-\mu}{\sqrt{|\mu|}}\,JE_j(f)\\
    &=&
    \frac{1}{\sqrt{|\mu|}}\,(\overline{Hess
}f(N,JE_j)+(1-\mu)JE_j(f))\\
 &=&\frac{1}{\sqrt{|\mu|}}\,((s-1)\langle
\overline{{\rm grad}}f,JE_j\rangle
+(1-\mu)JE_j(f))\\
&=&\frac{s-\mu}{\sqrt{|\mu|}}\,JE_j(f)=\frac{\mu-s}{\sqrt{|\mu|}}\,g_\mu(JC,E_j),
\end{eqnarray*}
and therefore $\nabla^{\mu}_{V^{\mu}}C=
\frac{\mu-s}{\sqrt{|\mu|}}\,JC$.
\end{proof}
\end{prop}

\begin{prop}\label{hessianos grado s}Let $V^\mu$ be the unit Hopf vector field in $(S^{2m+1},g_\mu)$ with $\mu\neq
0$ and $C={\rm grad}^\mu f-\varepsilon_\mu V^\mu(f)V^\mu$, where
$f$ is a harmonic and homogeneous polynomial of degree $s$
verifying (\ref{JHess}), then:
\begin{eqnarray*}
&&\textrm{a) }(Hess
E_{g_\lambda})_{V^\mu}(C)=\sqrt{\frac{|\lambda|}{|\mu|}}\int_{S^{2m+1}}\Big(\big(2ms((1-2m)\mu+\frac{(s-\mu)^2}{\lambda}
-4(s-1)^2\\
&&\hspace{2.5in}+2s)-4 s^2(s-1)^2)\big)f^2 +\|\overline{Hess}f\|^2\Big)\,{\rm dv}_\mu,\\
&&\textrm{b) }(Hess
E)_{V^\mu}(C)=\int_{S^{2m+1}}\Big(\big(2ms((2-2m)\mu+\frac{s^2}{\mu}
-4(s-1)^2)\\
&&\hspace{3in}-4 s^2(s-1)^2)\big)f^2
 +\|\overline{Hess}f\|^2\Big)\,{\rm dv}_\mu,\\
&&\textrm{c) }(Hess F)_{V^{\mu}}(C)=(1+|\mu|)^{m-2}\int_{S^{2m+1}}
    \Big(\big(2ms\ h(m,s,\lambda,\mu)-4s^2(s-1)^2\big)f^2\\
    &&\hspace{4in}+\|\overline{Hess}f\|^2\Big)\,{\rm
    dv}_\mu,
\end{eqnarray*}
where
$$h(m,s,\lambda,\mu)=\mu((2-2m)(1+|\mu|)+2\,\varepsilon_\mu\,(1+m-2s))
+\frac{s^2}{\mu}-4(s-1)^2+\varepsilon_\mu\,s^2.$$
\begin{proof}
By Proposition \ref{HesHopf}, to compute the Hessians of the
functionals when they act on the vector field $C$, we need to know
$\|\nabla^\mu C\|^2$, but
\begin{displaymath}
\|\nabla^\mu C\|^2=\sum_{i,j=1}^{2m}(B_i^j)^2 +
\mu\,\|C\|^2+\varepsilon_\mu\,\|\nabla^{\mu}_{V^{\mu}}C\|^2,
\end{displaymath}
where
\begin{eqnarray*}
    B_i^j&=&g_{\mu}(\nabla^{\mu}_{E_i}C,E_j)=(Hess )^{\mu}f(E_i,E_j)
    -\varepsilon_{\mu}\,g_{\mu}(\nabla_{E_i}^{\mu}(V^{\mu}(f)V^{\mu}),E_j)\\
    &=&\overline{Hess}f(E_i,E_j)-sf
    g(E_i,E_j)-V(f)g(JE_i,E_j).
\end{eqnarray*}
Therefore,
\begin{eqnarray*}
\sum_{i,j=1}^{2m}(B_i^j)^2&=&\sum_{i,j=1}^{2m}(\overline{Hess}f(E_i,E_j)-s f g(E_i,E_j)-V(f)g(JE_i,E_j))^2\\
&=&\sum_{i=1}^{2m}(\overline{Hess}f(E_i,E_i)-sf)^2+
\sum_{i=1}^m(\overline{Hess}f(E_i,JE_i)-V(f))^2\\
&&+\sum_{i=1}^m(\overline{Hess}f(JE_i,E_i)+V(f))^2+\sum_{{\rm default}}(\overline{Hess}f(E_i,E_j))^2\\
&=&\|\overline{Hess}f\|^2-2\sum_{i=1}^{2m}(\overline{Hess
}f(E_i,V))^2-2\sum_{i=1}^{2m}(\overline{Hess
f}(E_i,N))^2\\
&&-2(\overline{Hess}f(N,N))^2-2(\overline{Hess }f(N,V))^2+2m
V(f)^2+2ms^2f^2. \end{eqnarray*}

Now, since
$$\overline{Hess}f(E_i,V)=\overline{Hess}f(JE_i,N)=(s-1)JE_i(f),\quad
 \overline{Hess}f(E_i,N)=(s-1)E_i(f),$$
$$\overline{Hess}f(N,N)=s(s-1)f\quad \textrm{and}\quad
\overline{Hess}f(N,V)=(s-1)V(f),$$
 we have that,
\begin{displaymath}
 \sum_{i,j=1}^{2m}(B_i^j)^2=\|\overline{Hess
}f\|^2-4(s-1)^2\|C\|^2+2s^2f^2(m-(s-1)^2)+2(m-(s-1)^2)V(f)^2
\end{displaymath}
and
\begin{displaymath}
\|\nabla^\mu C\|^2=\|\overline{Hess
}f\|^2+(\mu+\frac{(s-\mu)^2}{\mu}-4(s-1)^2)\|C\|^2+2(m-(s-1)^2)(s^2f^2+V(f)^2).
\end{displaymath}
Moreover, since $V$ is a Killing vector field,
\begin{eqnarray*}
\int_{S^{2m+1}}(V(f))^2{\rm
dv}_\mu&=&-\int_{S^{2m+1}}V(V(f))\,{\rm
dv}_\mu=s^2\int_{S^{2m+1}}f^2\,{\rm dv}_\mu,\\
\int_{S^{2m+1}}\|C\|^2\,{\rm
dv}_\mu&=&(2ms+\frac{s^2}{\mu})\int_{S^{2m+1}}f^2\,{\rm
dv}_\mu-\frac{1}{\mu}\int_{S^{2m+1}}(V(f))^2\,{\rm
dv}_\mu\\
&=&2ms\int_{S^{2m+1}}f^2\,{\rm dv}_\mu,
\end{eqnarray*}
from where the result holds.

\end{proof}
\end{prop}

Let $(x_1,\ldots,x_m,x_{m+1},\ldots,x_{2m+2})$ be coordinates in
 $\mathbb{R}^{2m+2}$. If we denote $z_j=x_j+ i\, x_{m+1+j}$, then
$(z_1,\ldots,z_{m})\in \mathbb{C}^m$ and the complex structure of
$\mathbb{R}^{2m+2}$ is given by
$J(\partial_{x_j})=\partial_{x_{m+1+j}}$,
$J(\partial_{x_{m+1+j}})=-\partial_{x_j}$. We are going to compute
the Hessians of the functionals in the direction of vector fields
$C$ such that the polynomial $f$ depends only on two variables
$x_j$, $x_{m+1+j}$ that we will represent by $x, y$.\par It is
easy to see that the polynomials of even degree
$f=\sum_{i=0}^{s}(-1)^i \left(\begin{array}{c}
2s\\2i\end{array}\right)x^{2(s-i)}y^{2i}$ verifies  the hypothesis
of the above Proposition. Then,
we have to compute
\begin{displaymath}
\int_{S^{2m+1}}f^2\,{\rm dv}_\mu \qquad \textrm{ and } \qquad
\int_{S^{2m+1}}\|\overline{Hess}f\|^2\,{\rm dv}_\mu.
\end{displaymath}

In order to do so, since
$$\left(\begin{array}{c} 2s\\2i\end{array}\right)=\left(\begin{array}{c}
2s-1\\2i\end{array}\right)+\left(\begin{array}{c}
2s-1\\2i-1\end{array}\right)$$ if we take
\begin{displaymath}
Y=\sum_{i=0}^{s-1}(-1)^i \left(\begin{array}{c}
2s-1\\2i\end{array}\right)x^{2(s-i)-1}y^{2i}\partial_x +
\sum_{i=1}^s(-1)^i \left(\begin{array}{c}
2s-1\\2i-1\end{array}\right)x^{2(s-i)}y^{2i-1}\partial_y,
\end{displaymath}
then $f^2=\langle fY,N\rangle$ and
\begin{displaymath}
\int_{S^{2m+1}}f^2\,{\rm
dv}_\mu=\sqrt{|\mu|}\int_{B^{2m+2}}\overline{{\rm
div}}(fY)\,\omega_{2m+2}=\sqrt{|\mu|}\int_{B^{2m+2}}Y(f)\,\omega_{2m+2},
\end{displaymath}
 since $\overline{{\rm
div}}(Y)=0$. Here $\omega_{2m+2}$ is the volume element on $\mathbb{R}^{2m+2}$.\\
In addition, it is easy to see applying induction that
$Y(f)=2s(x^2+y^2)^{2s-1}$ and then,
\begin{displaymath}
\int_{S^{2m+1}}f^2\,{\rm
dv}_\mu=2s\sqrt{|\mu|}\int_{B^{2m+2}}(x^2+y^2)^{2s-1}\,\omega_{2m+2}.
\end{displaymath}
On the other hand, easy computations show that
$\|\overline{Hess}f\|^2= 8s^2(2s-1)^2 (x^2+y^2)^{2s-2}$, and
taking
$\tilde{Y}=8s^2(2s-1)^2(x(x^2+y^2)^{2s-3}\partial_x+y(x^2+y^2)^{2s-3}\partial_y)$,
\begin{eqnarray*}
\int_{S^{2m+1}}\|\overline{Hess}f\|^2\,{\rm dv}_\mu
&=&\sqrt{|\mu|}\int_{B^{2m+2}}\overline{{\rm
div}}(\tilde{Y})\,\omega_{2m+2}\\
&=& 8s^2(2s-1)^2
(4s-4)\sqrt{|\mu|}\int_{B^{2m+2}}(x^2+y^2)^{2s-3}\,\omega_{2m+2}.
\end{eqnarray*}

Moreover, it is not difficult to see that
\begin{displaymath}
\int_{B^{2m+2}}(x^2+y^2)^{2s-1}\,\omega_{2m+2}=\frac{(2s-1)(2s-2)}{(m+2s-1)(m+2s)}\int_{B^{2m+2}}(x^2+y^2)^{2s-3}\,\omega_{2m+2},
\end{displaymath}
and then,
\begin{displaymath}
\int_{S^{2m+1}}\|\overline{Hess}f\|^2\,{\rm dv}_\mu=
8s(2s-1)(m+2s-1)(m+2s)\int_{S^{2m+1}}f^2\,{\rm dv}_\mu.
\end{displaymath}
As a consequence,
\begin{lem}\label{hessianospolinomios2s} Let $V^\mu$ be the unit Hopf vector field on $(S^{2m+1},g_\mu)$
with $\mu\neq 0$. Then  for each $s>0$ there exists a vector field
$C_{2s}={\rm grad}^\mu f_{2s}-\varepsilon_\mu
V^{\mu}(f_{2s})V^\mu$ orthogonal to $V$ such that
 \begin{eqnarray*}
&&\textrm{a) }(Hess
E_{g_\lambda})_{V^\mu}(C_{2s})=\sqrt{\frac{|\lambda|}{|\mu|}}\,e_\lambda(\mu,m,s)
\int_{S^{2m+1}}\|C_{2s}\|^2\,{\rm dv}_\mu.\\
&&\textrm{b) }(Hess
E)_{V^\mu}(C_{2s})=\frac{2}{\mu}(\mu^2(1-m)+\mu(2s-1)(m+1)+2s^2)\int_{S^{2m+1}}\|C_{2s}\|^2\,{\rm
dv}_\mu.\\
&&\textrm{c)
}(HessF)_{V^\mu}(C_{2s})=\frac{2}{\mu}(1+|\mu|)^{m-2}f(s,m,\mu)\int_{S^{2m+1}}\|C_{2s}\|^2\,{\rm
 dv}_\mu.
 \end{eqnarray*}
Where
\begin{eqnarray*}
&&e_\lambda(\mu,s,m)=\mu(1-2m)+\frac{(2s-\mu)^2}{\lambda}+2(2s-1)(m+1)+4s,\\
&&f(s,m,\mu)=\mu^2(1-m)(1+|\mu|)+\mu|\mu|(1+m-4s)+\mu((2s-1)(m+1)+2\,\varepsilon_\mu\,s^2)
 +2s^2.
\end{eqnarray*}
 \end{lem}
 Using these expressions, we can show that
 \begin{thm}\label{inest todo Lorentz} On $(S^{2m+1},g_\mu)$  the unit Hopf vector
 fields are unstable critical points of the energy and the volume for all
 $\mu<0$, and they are unstable as a critical points of the generalized
energies $E_{g_\lambda}$
 for all $\lambda <0$ and $\mu\neq 0$.
 \begin{proof}
 Using b) of Lemma \ref{hessianospolinomios2s} for each
 $s>0$ there exists a vector field $C_{2s}$ such that
\begin{displaymath}
(Hess
E)_{V^\mu}(C_{2s})=\frac{2}{\mu}(\mu^2(1-m)+\mu(2s-1)(m+1)+2s^2)\int_{S^{2m+1}}\|C_{2s}\|^2\,{\rm
dv}_\mu.
\end{displaymath}
For each $\mu<0$ fixed,
$\frac{2}{\mu}(\mu^2(1-m)+\mu(2s-1)(m+1)+2s^2)$ goes to $-\infty$
as $s$ grows, so there exists $s>0$ such that $(Hess
E)_{V^\mu}(C_{2s})<0$ and therefore $V^\mu$ is unstable for the
energy.\par Analogously, we obtain the instability with respect to
the other functionals.
\end{proof}
\begin{nota} For $\mu<0$ and $\lambda>0$, using a) of Proposition
\ref{HesHopf} we have that $(Hess E_{g_\lambda})_{V^\mu}(A)\geq 0$
for all $A \in V^\bot$ and then $V^\mu$ is stable for the
generalized energy $E_{g_\lambda}$.
\end{nota}

\end{thm}
\section{Riemannian Berger spheres}
In \cite{GMH}, the authors studied the stability of Hopf vector
fields in Berger Riemannian spheres with respect to the energy,
the volume and the generalized energy. This problem is completely
solve for the energy and the volume, but as concerns the
generalized energy, there exist values of $\lambda$ and $\mu$ for
which the behavior of Hopf vector fields is unknown. For the
$3$-sphere, it is shown that
\begin{prop}[\cite{AH}] On $(S^3,g_\mu)$ with $\mu>0$, if
$\lambda>(\mu-2)^2/\mu$ then the Hopf vector field $V^\mu$ is an
unstable critical point of $E_{g_\lambda}$.
\end{prop}
\begin{thm}[\cite{AH}] On $(S^3,g_\mu)$ with $\mu>0$, the Hopf
vector field $V^\mu$ is stable as a critical point of the
functionals $E_{g_\lambda}$ in the followings cases:
\begin{itemize}
\item[a) ]If $\mu\leq 8/3$, for $\lambda\leq (\mu-2)^2/\mu$,
\item[b) ]If $8/3<\mu\leq 4$, for $\lambda\leq (\mu-3)^2/(\mu-2)$,
\item[c) ]If $\mu>4$, for $\lambda\leq \mu-4$.
\end{itemize}
\end{thm}
As we have seen in the previous section, the simultaneous
eigenfunctions of the Laplacian and of the vertical Laplacian
 have allow us to construct directions such
that the Hessians take negative values when they act on them. We
are going to use the same idea, but now, using the special
structure of the $3$-sphere. We construct vector fields $A=a_1
E_1+a_2 E_2$, where if $i$, $j$, $k$, represent the imaginary unit
quaternions then $\{V^\mu, E_1=jN,E_2=kN\}$ is an adapted
$g^\mu$-orthonormal frame and where $a_1$, $a_2$ are
eigenfunctions of the Laplacian.\par

\begin{figure}[hbt]
\centering
\includegraphics{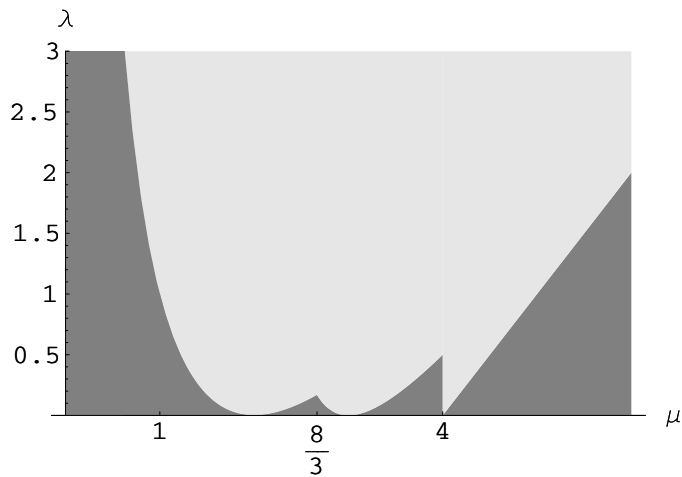}
\caption{The light gray region is the subset of
$\mathbb{R}^+\times \mathbb{R}^+$ of pairs $(\mu,\lambda)$ such
that $V^\mu$ is unstable as a critical point of $E_{g_\lambda}$.
The stability domain is painted in dark gray.} \label{grafica}
\end{figure}

\begin{prop}\label{Inestdim3-2} On $(S^3,g_\mu)$ with $\mu>0$, if
 $\lambda>(\mu-3)^2/(\mu-2)$ and $\mu>2$, the Hopf vector filed $V^\mu$ is unstable as a critical point of
  $E_{g_\lambda}$.
\begin{proof}
We take $A=a_1\, E_1 +a_2\, E_2$ with $a_2=V(a_1)$ and $a_1$ an
eigenfunction associated to the first eigenvalue of the Laplacian
$\lambda_1=3$ satisfying that $V(V(a_1))=-a_1$. That is to say,
$a_1,a_2 \in \mathcal{P}_1^1$.\par If we compute $\nabla_V A$ we
obtain that,
\begin{eqnarray*}
\nabla_V A&=&\nabla_V (a_1 E_1+a_2
E_2)=V(a_1)E_1+a_1(-E_2)+V(a_2)E_2+a_2 E_1\\
&=&(V(a_1)+a_2)E_1+(V(a_2)-a_1)E_2=2 a_2E_1-2a_1E_2=-2 JA,
\end{eqnarray*}
and therefore
\begin{displaymath}
\nabla^\mu_{V^\mu}A=\frac{1}{\sqrt{\mu}}(\nabla_V
A+(\mu-1)\nabla_A V)=\frac{\mu-3}{\sqrt{\mu}}JA.
\end{displaymath}
By Proposition 4.4 of \cite{GMH}
\begin{eqnarray*}
(Hess
E_{g_{\lambda}})_{V^{\mu}}(A)&=&\sqrt{\lambda/\mu}\int_{S^{3}}
   \Big((\mu-4-\lambda+\frac{(\mu-3-\lambda)^2}{\lambda})\|A\|^2+\frac{1}{2}\|\bar D^C
    A\|^2_{V^\bot}\Big)\,{\rm dv}_\mu\\
    &=&\sqrt{\lambda/\mu}\int_{S^{3}}\Big((2-\mu+\frac{(\mu-3)^2}{\lambda})\|A\|^2+\frac{1}{2}\|\bar D^C
    A\|^2_{V^\bot}\Big)\,{\rm dv}_\mu,
\end{eqnarray*}
where
\begin{eqnarray*}
\int_{S^3}\frac{1}{2}\|\bar D^C
    A\|^2_{V^\bot}\,{\rm
    dv}&=&\int_{S^3}\Big((B_1^1)^2+(B_2^1)^2+(B_1^2)^2+(B_2^2)^2\Big)\,{\rm
    dv}\\
    &&-2\int_{S^3}(B_2^1B_1^2-B_2^2 B_1^1)\,{\rm dv},
\end{eqnarray*}
and
\begin{eqnarray*}
&& B_1^1=E_1(a_1), \qquad B_1^2=E_1(a_2),\\
&& B_2^1=E_2(a_1), \qquad B_2^2=E_2(a_2).
\end{eqnarray*}
So, it is enough to prove that
\begin{displaymath}
\int_{S^3}\frac{1}{2}\|\bar D^C
    A\|^2_{V^\bot}\,{\rm dv}_\mu=\sqrt{\mu}\int_{S^3}\frac{1}{2}\|\bar D^C
    A\|^2_{V^\bot}\,{\rm dv}=0.
\end{displaymath}

Since $E_1$ and $E_2$ are Killing vector fields
\begin{eqnarray*}
\int_{S^{3}}\sum_{i,j=1}^2(B_i^j)^2\,{\rm
dv}&=&\int_{S^{3}}\Big((E_1(a_1))^2+(E_1(a_2))^2+(E_2(a_1))^2+(E_2(a_2))^2\Big)\,{\rm
dv}\\
&=&-\int_{S^3}\Big(a_1 E_1(E_1(a_1))+a_2 E_1(E_1(a_2))\\
&&\hspace{1in}+a_1
E_2(E_2(a_1))+a_2 E_2(E_2(a_2))\Big)\,{\rm dv}\\
&=&\int_{S^3} \sum_{i=1}^2 a_i(\Delta (a_i)+V(V(a_i)))\,{\rm dv}=
2\int_{S^3}\|A\|^2\,{\rm dv}.
\end{eqnarray*}
On the other hand,
\begin{eqnarray*}
\int_{S^3}(B_2^1B_1^2-B_2^2 B_1^1)\,{\rm
dv}&=&\int_{S^3}(E_2(a_1)E_1(a_2)-E_2(a_2) E_1(a_1))\,{\rm dv}\\
&=&-\int_{S^3}(a_2 E_1(E_2(a_1))-a_2 E_2(E_1(a_1)))\,{\rm dv}\\
&=&-\int_{S^3}a_2 [E_1,E_2](a_1)\,{\rm dv}=2\int_{S^3}a_2
V(a_1)\,{\rm dv}\\
&=&2\int_{S^3}a_2^2\,{\rm dv}=\int_{S^3}(a_1^2+a_2^2)\,{\rm dv},
\end{eqnarray*}
since
\begin{displaymath}
\int_{S^3}a_2 V(a_1)\,{\rm dv}=-\int_{S^3}a_1 V(a_2)\,{\rm
dv}=\int_{S^3}a_1^2\,{\rm dv}.
\end{displaymath}
Consequently,
\begin{displaymath}
\int_{S^3}\frac{1}{2}\|\bar D^C
    A\|^2_{V^\bot}\,{\rm
    dv}=2\int_{S^3}\|A\|^2\,{\rm dv}-2\int_{S^3}\|A\|^2\,{\rm dv}=0.
\end{displaymath}
\end{proof}
\end{prop}
Using similar arguments, we can prove that
\begin{prop}\label{Inestdim3-3} On $(S^3,g_\mu)$ with $\mu>0$, if
$4-\mu+(\mu-4)^2/\lambda<0$, or equivalently, if $\lambda>\mu-4$
and $\mu>4$, then  $V^\mu$ is an unstable critical point of
$E_{g_\lambda}$.
\begin{proof}
Now, we take $A=a_1\, E_1 +a_2\, E_2$ with $a_2=V(a_1)/2$ and
$a_1$ an eigenfunction associated to the eigenvalue of the
Laplacian $\lambda_2=8$ verifying that $V(V(a_1))=-4 a_1$. That is
to say, $a_1,a_2 \in \mathcal{P}^2_2$.\par If we compute $\nabla_V
A$ we obtain,
\begin{eqnarray*}
\nabla_V A&=&3 a_2E_1-3a_1E_2=-3 JA,
\end{eqnarray*}
and then
\begin{displaymath}
\nabla^\mu_{V^\mu}A=\frac{\mu-4}{\sqrt{\mu}}JA.
\end{displaymath}
By Proposition 4.4 of \cite{GMH}
\begin{eqnarray*}
(Hess E_{g_{\lambda}})_{V^{\mu}}(A)=&=&
\sqrt{\lambda/\mu}\int_{S^{3}}\Big((4-\mu+\frac{(\mu-4)^2}{\lambda})\|A\|^2+\frac{1}{2}\|\bar
D^C
    A\|^2_{V^\bot}\Big)\,{\rm dv}_\mu.
\end{eqnarray*}
Using the same arguments that in the above Proposition
\begin{displaymath}
\int_{S^{3}}\sum_{i,j=1}^2(B_i^j)^2\,{\rm
dv}=4\int_{S^3}(a_1^2+a_2^2)\,{\rm dv}
\end{displaymath}
and
\begin{displaymath}
\int_{S^3}(B_2^1B_1^2-B_2^2 B_1^1)\,{\rm
dv}=4\int_{S^3}a_2^2\,{\rm dv}=2\int_{S^3}(a_1^2+a_2^2)\,{\rm dv}.
\end{displaymath}
Therefore,
\begin{displaymath}
\int_{S^3}\frac{1}{2}\|\bar D^C
    A\|^2_{V^\bot}\,{\rm
    dv}=0.
\end{displaymath}
\end{proof}
\end{prop}
These results jointly with those shown in $\cite{AH}$, solve
completely the problem of determining the behaviour of Hopf vector
fields in Berger Riemannian $3$-spheres. These results can be
represented graphically in $\mathbb{R}^+\times \mathbb{R}^+$, as
can be seen in Figure \ref{grafica}.\par

For spheres of upper dimension, we can use the vector fields
$C_{2s}$ introduced in the previous section and we can also,
following the idea used to solve the problem in dimension $3$,
construct vector fields $A=a_1 A_a +a_2 A_{Ja}$, with $a_1$ and
$a_2$ simultaneous eigenfunctions of the Laplacian and of the
vertical Laplacian. With these directions we can improve the
results stayed in $\cite{GMH}$, but unfortunately, it is not
sufficient to solve the problem.


\end{document}